\documentclass[12pt]{article}

\usepackage{amsfonts}

\newcommand{\sect}[1]{\setcounter{equation}{0}\section{#1}}
\newcommand{\subsect}[1]{\subsection{#1}}

\def\be{\begin{equation}}
\def\ee{\end{equation}}
\def\bea{\begin{eqnarray}}
\def\eea{\end{eqnarray}}

\def\1{\'{\i}}                           
\def\Re{{\mathbb R}}

\def\jp{J_+}
\def\jm{J_-}
\def\jj{J_3}

\def\para{\omega^2}

\def\kk{K}

\def\s{{\sigma}}

\def\co{\Delta}

\def\J{{J}} 
 
\def\s{{\sigma}} 

\def\>#1{{\bf #1}}                 
\def\1{\'{\i}}                           
\def\R{\rm I\kern-.2em R} 
\def\C{\rm I\kern-.5em C}

\def\tfrac#1#2{ {\scriptstyle { \frac {#1}{#2}}}}         
\def\pois#1#2{\left\{ {#1},{#2} \right\}}

\begin{document}
 
\thispagestyle{empty}
\

\begin{center} {\LARGE{\bf{Integrable deformations of Hamiltonian
}}}


{\LARGE{\bf{systems and $q$-symmetries
\footnote{Contribution to the III International Workshop on {Classical
and Quantum Integrable Systems}. Edited by L.G. Mardoyan, G.S.
Pogosyan  and A.N. Sissakian. Joint Institute for Nuclear
Research, Dubna, pp. 15--25, (1998).}
}}} 

\vspace{1.4cm} 

ANGEL BALLESTEROS and FRANCISCO J. HERRANZ\vspace{.2cm}\\
{\it Departamento de F\1sica, Universidad de Burgos\\
Pza. Misael Ba\~nuelos s.n., 09001-Burgos, Spain} 
\vspace{.3cm}

\end{center} 
  
\bigskip

\begin{abstract} 

The complete integrability of the hyperbolic Gaudin Hamiltonian and
other related integrable systems is shown to be easily derived by taking
into account their $sl(2,\Re)$ coalgebra symmetry. By using the
properties induced by such a coalgebra structure, it can be proven that
the introduction of any quantum deformation of the $sl(2,\Re)$ algebra
will provide an integrable deformation for such systems. In particular,
the Gaudin Hamiltonian arising from the non-standard quantum deformation
of the $sl(2,\Re)$ Poisson algebra is presented, including the explicit
expressions for its integrals of motion. A completely integrable system of
nonlinearly coupled oscillators derived from this deformation is also
introduced.

\end{abstract}


\sect{Introduction: on $q$-Poisson coalgebras}

Let us consider the $sl(2)$ Poisson-Lie algebra 
\be
\pois{J_3}{J_\pm}=\pm\,2\,J_\pm, \qquad 
\pois{J_-}{J_+}=\alpha\,J_3. \qquad 
\label{hg}
\ee
(where $\alpha$ is a real positive parameter and the usual $\alpha=1$
case can be easily get by performing the change of basis $J_\pm\to
J_\pm/\sqrt{\alpha}$). The Casimir function for this algebra reads
 \be
C(J_3,J_\pm)=\frac{\alpha}{4} J_3^2 - \,J_+J_-.
\label{hgb}  
\ee
The Poisson algebra (\ref{hg}) is endowed with a coalgebra structure
\cite{Hopf,Sweed} by means of the ``primitive" coproduct 
$\co:sl(2) \rightarrow sl(2)\otimes sl(2)$ defined as follows:
\be
\co(J_3)=J_3\otimes 1 + 1\otimes J_3 ,\qquad
\co(J_\pm)= J_\pm\otimes 1 + 1\otimes J_\pm.
\label{prcop}
\ee
Compatibility between (\ref{prcop}) and (\ref{hg}) means that $\co$ is a
Poisson map: the three functions defined through (\ref{prcop}) close also
a Poisson $sl(2)$ algebra. The map (\ref{prcop}) can be introduced in a
similar manner for any Lie algebra. The coalgebra structure,
together with the known relevance of Lie algebras in the context of
Liouville integrability \cite{OP,Per}, are the keystones of the
universal construction of integrable systems introduced in \cite{BR} that
we shall summarize in the next Section.

From that point of view, the search for $q$-Poisson algebras can
be seen as the problem of finding deformations of a given algebra endowed
with compatible deformations of the coproduct (\ref{prcop}) preserving it
as a Poisson map. In fact, two relevant and distinct structures
deforming $sl(2,\Re)$ appeared in quantum group literature during last
years (see, for instance, \cite{CP}), and can be realized as $q$-Poisson
algebras as follows:

\noindent a) The {\it ``standard"} deformation \cite{Dr,Ji}
$sl_q(2)$, which is given by the following deformed Poisson brackets
$(q=e^z)$:
\be
\pois{\J_3}{\J_+}=2\,\J_+,\quad \pois{\J_3}{\J_-}=-2\,\J_-,\quad 
\pois{\J_-}{\J_+}=\frac{\sinh z\J_3}{z},
\label{lc} 
\ee
which are compatible with the deformed coproduct
\bea
&& \co(\J_3) =\J_3\otimes 1 + 1 \otimes \J_3  ,\cr  
&&  \co(\J_+) =\J_+\otimes
e^{\tfrac{z}{2}\J_3} + e^{-\tfrac{z}{2}\J_3} \otimes \J_+ , \label{lb} \\
&&  \co(\J_-) =\J_-\otimes
e^{\tfrac{z}{2}\J_3} + e^{-\tfrac{z}{2}\J_3} \otimes \J_- ,\nonumber
\eea
in the sense that that (\ref{lb}) is a Poisson algebra
homomorphism with respect to (\ref{lc}). The function
\be
{\cal C}_z^{(s)}(\J_3,\J_\pm)=\left(\frac{\sinh (\tfrac{z}{2} \J_3)}{z}
\right)^2 - \J_+\,\J_-, \label{le}
\ee 
is the deformed Casimir for this Poisson coalgebra.

\noindent b) The {\it ``non-standard"} quantum
$sl(2,\Re)$ deformation \cite{Demidov,Ohn} is defined by:
\bea
&&\Delta(\jm)=  \jm \otimes 1+
1\otimes \jm ,\cr
&&\Delta(\jp)=\jp \otimes e^{z \jm} + e^{-z \jm} \otimes
\jp, \cr
&&\Delta(\jj)=\jj \otimes e^{z \jm} + e^{-z \jm} \otimes \jj ,
\label{ga}
\eea
\be 
\{\jj,\jp\}=2 \jp \cosh z\jm, \qquad 
 \{\jj,\jm\}=-2\frac {\sinh z\jm}{z},\qquad
\{\jm,\jp\}=4 \jj . 
\label{gb}
\ee 
The deformed Casimir reads
\be
{\cal C}_z^{(n)}(\J_3,\J_\pm)=\jj^2  - \frac {\sinh z\jm}{z} \jp.
\label{gc}
\ee

The aim of this contribution is to present an application of this latter
$q$-Poisson algebra in the construction of an integrable deformation of
the Gaudin Hamiltonian. This deformation can be also interpreted in the
context of a chain of nonlinearly coupled oscillators or, equivalently,
in relation with a certain integrable perturbation of the motion of a 
particle under any central
potential in the $N$-dimensional Euclidean space. Through these examples
we will show an intrinsic connection between quantum deformations and
nonlinear interactions depending on the momenta. Results concerning the
standard deformation of the Gaudin-Calogero system \cite{Cal} can be found
in \cite{BR,BCR}.


\sect{The formalism}

By following \cite{BR}, we can state the following result: any
coalgebra $(A,\co)$ with Casimir element $C$ can be considered as the
generating symmetry of a large family of integrable
systems. We shall consider here classical mechanical systems only and,
consequently, we shall make use of Poisson realizations $D$ of the
algebra $A$ of the form $D:A\rightarrow C^\infty (q,p)$. However, we
recall that the formalism is also directly applicable to quantum
mechanical systems. The constructive procedure is as follows.

Let
$(A,\Delta)$ be a (Poisson) coalgebra with generators $X_i$
$(i=1,\dots,l)$ and Casimir function $C(X_1,\dots,X_l)$. This means that
the coproduct $\Delta:A\rightarrow A\otimes A$ is a Poisson map. Let us
consider the $N$-th coproduct map $\co^{(N)}$ 
\be
\co^{(N)}:A\rightarrow A\otimes A\otimes \dots^{N)}\otimes A ,
\ee
which is obtained (see \cite{BR}) by applying recursively the
two-coproduct
$\co^{(2)}\equiv \co$ in the form
\be
\co^{(N)}:=(id\otimes id\otimes\dots^{N-2)}\otimes id\otimes
\co^{(2)})\circ\co^{(N-1)}.
\label{fl}
\ee
 By taking into account that the $m$-th
coproduct $(m\leq N)$ of the Casimir $\co^{(m)}(C)$ can be embedded into
the tensor product of $N$ copies of $A$ as
\be
\co^{(m)}:A\rightarrow \{A\otimes A\otimes \dots^{m)}\otimes A\} \otimes
\{1 \otimes 1\otimes \dots^{N-m)}\otimes 1\},
\ee
it can be shown that, \be
\pois{\co^{(m)}(C)}{\co^{(N)}(X_i)}_{A\otimes
A\otimes\dots^{N)}\otimes A}=0, \qquad
i=1,\dots,l, \qquad 1\leq m\leq N .
\label{za}
\ee

With this in mind it can be proven \cite{BR} that, if ${\cal
H}$ is an {\it arbitrary} (smooth) function of the
generators of $A$, the $N$-particle Hamiltonian defined on $A\otimes
A\otimes\dots^{N)}\otimes A$ as the $N$-th coproduct of ${\cal H}$
\be
H^{(N)}:=\co^{(N)}({\cal{H}}(X_1,\dots,X_l))=
{\cal{H}}(\co^{(N)}(X_1),\dots,\co^{(N)}(X_l)),
\label{htotg}
\ee
fulfils
\be
\pois{C^{(m)}}{H^{(N)}}_{A\otimes
A\otimes\dots^{N)}\otimes A}=0, \qquad 1\leq m\leq N,
\label{za1}
\ee
where the $N$ functions $C^{(m)}$ ($m=1,\dots,N$)
are defined through the $m$-th coproducts of the Casimir $C$
\be
C^{(m)}:= \co^{(m)}(C(X_1,\dots,X_l))=
C(\co^{(m)}(X_1),\dots,\co^{(m)}(X_l)) ,
\label{Ctotg}
\ee
and all the integrals of motion $C^{(m)}$ are in involution
\be
\pois{C^{(m)}}{C^{(n)}}=0, \qquad \forall\,m,n=1,\dots,N.
\label{cor3}
\ee

Therefore, provided a non-trivial realization of $A$ on a one-particle
phase space is given, the
$N$-particle Hamiltonian $H^{(N)}$ will be a function of $N$ canonical
pairs $(q_i,p_i)$ and is, by construction, completely  integrable with
respect to the ordinary Poisson bracket 
\be
\{f,g\}=\sum_{i=1}^N\left(\frac{\partial f}{\partial q_i}
\frac{\partial g}{\partial p_i}
-\frac{\partial g}{\partial q_i} 
\frac{\partial f}{\partial p_i}\right) .
\label{poisbra}
\ee
Moreover, its
constants of motion will be given by the $C^{(m)}$ functions, all of
them functionally independent since each of them depends on the first $m$
pairs $(q_i,p_i)$ of canonical coordinates. Note that with such
one-particle realizations the first Casimir $C^{(1)}$ will be a number,
and we are left with $N-1$ constants of motion with respect to $H^{(N)}$.

Let us stress now that the previous construction is valid for
quantum algebras with no extra assumptions. Quantum algebras are also
(deformed) coalgebras $(A_z,\Delta_z)$, and any function of the generators
of a given quantum algebra with Casimir element $C_z$ will provide, under
a chosen deformed representation, a completely integrable Hamiltonian.
Therefore, the obtention of quantum algebras is a direct method to get
integrable deformations of those Hamiltonian systems with underlying
coalgebra symmetry.


\sect{Oscillator chains from $sl(2,\Re)$ coalgebras}

It is well-known that $sl(2,\Re)$ can be considered as a dynamical
algebra for the one-dimensional harmonic oscillator. As we shall see in
the sequel, the $sl(2,\Re)$ coalgebra will give us the complete set of
integrals of the motion of a chain of $N$ independent harmonic
oscillators with the same frecuency $\omega$.
When the non-standard quantum
deformation is considered, a new integrable oscillator
chain with long range interactions depending on the momenta is obtained.


\subsect{The non-deformed case}

Let us consider the $sl(2,\Re)$ coalgebra (\ref{hg}) and (\ref{hgb}) with
$\alpha=4$. A one-particle phase space realization of this Poisson
algebra with vanishing Casimir is given by the functions:
\be 
 f_-^{(1)}=D(\jm)=q_1^2, \qquad 
f_+^{(1)}=D(\jp)= p_1^2, \qquad
 f_3^{(1)}=D(\jj)=q_1 p_1 .
\label{fd}
\ee 
From it, the harmonic oscillator Hamiltonian is
recovered if the following linear function of the generators of
$sl(2,\Re)$
\be
{\cal H}= \jp+\para  \jm ,
\label{hsl}
\ee
is represented through (\ref{fd}):
\be
H^{(1)}=D({\cal H})=   p_1^2 +\para q_1^2 .
\label{fe}
\ee

Now, the representation $(D\otimes D)$ when applied onto the
primitive coproduct (\ref{prcop}), leads to the following two-particle phase
space realization of $sl(2,\Re)$ \be 
 f_-^{(2)}=q_1^2+q_2^2, \qquad
 f_+^{(2)}= p_1^2 +  p_2^2, \qquad
 f_3^{(2)}=q_1 p_1 + q_2 p_2 . 
\label{ff}
\ee
that, in turn, gives rise to the uncoupled oscillator Hamiltonian:
\be
H^{(2)}= (D\otimes D)(\co^{(2)}({\cal H}))= f_+^{(2)} +\para  f_-^{(2)} 
=  p_1^2+  p_2^2+\para (q_1^2+q_2^2).
\label{fg}
\ee
Note that now the frecuency of both oscillators is the same. The phase
space realization of the coproduct of the Casimir will give us the
corresponding integral of motion \be
C^{(2)}=(D\otimes D)(\co^{(2)}({C}))=-({q_1}{p_2} -{q_2} {p_1})^2 ,
\label{fh}
\ee
that turns out to be the square of the angular momentum, as expected.

The construction of the $m$-dimensional functions from the $m$-th
coproduct is straightforwardly obtained by the induction (\ref{fl}):
\be f_-^{(m)}=
\sum_{i=1}^m q_i^2, \qquad
f_+^{(m)}
=\sum_{i=1}^m   p_i^2, \qquad
f_3^{(m)}
=\sum_{i=1}^m q_i p_i .
\label{fi}
\ee
From it, the uncoupled chain of $N$ harmonic oscillators (all of them with
the same frecuency) is obtained as the representation of the same
dynamical Hamiltonian (\ref{hsl}):
\be 
 H^{(N)}
=    f_+^{(N)} +\para  f_-^{(N)} 
=\sum_{i=1}^N(  p_i^2 +\para q_i^2)
\label{fj}
\ee
together with the integrals of motion $(m=2,\dots,N)$, that are deduced
from the $m$-th coproducts of the Casimir and are shown to be
\be
C^{(m)} =-\sum_{i<j}^m ({q_i}{p_j} - {q_j}{p_i})^2 .
\label{fk}
\ee
Note that the integrals $C^{(m)}$ given by the $sl(2,\Re)$  coalgebra are
just the quadratic Casimirs of the $so(m)$ algebras with $m=2,\dots,N$.
It is also  well known that the Hamiltonian (\ref{fj}) is $so(N)$
invariant, since it can be interpreted as the one for a particle moving
on the $N$-dimensional Euclidean space under the potential $\para\,r^2$.


\subsect{An integrable deformation from $U_z(sl(2,\Re))$}

Now, we introduce a  one-particle deformed phase space realization
of (\ref{gb}):
\bea 
&& f_-^{(1)}=D_z(\jm)=q_1^2, \qquad \cr
&& f_+^{(1)}=D_z(\jp)=\frac {\sinh z q_1^2}{z q_1^2}  p_1^2,   \label{gd}\\
&& f_3^{(1)}=D_z(\jj)=\frac {\sinh z q_1^2}{z q_1}  p_1 .
\nonumber
\eea 
This realization is also characterized by a vanishing deformed Casimir
function (\ref{gc})

Let us consider again the dynamical generator
${\cal H}= \jp +\para  \jm$. Under (\ref{gd}), we obtain the Hamiltonian
\be
H^{(1)}_z=D_z({\cal H})= \frac {\sinh z q_1^2}{z q_1^2} 
p_1^2  +\para q_1^2 .
\label{ge}
\ee

The corresponding
two-particle phase space realization of $U_z(sl(2,\Re))$ is obtained
from  both the coproduct (\ref{ga}) and (\ref{gd}):
\bea
&&f_-^{(2)}=q_1^2+q_2^2, \cr
&&f_+^{(2)}=
\frac {\sinh z q_1^2}{z q_1^2}  p_1^2  e^{z q_2^2} +
\frac {\sinh z q_2^2}{z q_2^2}  p_2^2  e^{-z q_1^2}, \cr
&&f_3^{(2)}=
\frac {\sinh z q_1^2}{z q_1 }  p_1   e^{z q_2^2} +
\frac {\sinh z q_2^2}{z q_2 }  p_2  e^{-z q_1^2} .
\label{gf}
\eea
The associated two-particle Hamiltonian is
\be 
 H^{(2)}_z=  \frac {\sinh z q_1^2}{z q_1^2}  p_1^2  e^{z q_2^2} +
\frac {\sinh z q_2^2}{z q_2^2}  p_2^2  e^{-z q_1^2} +\para
  ( q_1^2+q_2^2  ) ,
\label{gg}
\ee  
and the deformed coproduct for the deformed Casimir (\ref{gc}) reads 
\be
  C^{(2)}_z=-\frac {\sinh z q_1^2 \sinh z q_2^2}{z^2 q_1^2 q_2^2}
\left({q_1}{p_2} - {q_2}{p_1}\right)^2 e^{-z q_1^2}
e^{z q_2^2}.
\ee
If we rewrite (\ref{gg}) in the form
\bea
&& H^{(2)}_z= H^{(2)} + p_1^2 \left(\frac {\sinh z q_1^2}
{z q_1^2}e^{z q_2^2} -1 \right) + 
p_2^2 \left(\frac {\sinh z q_2^2}{z q_2^2}  e^{-z q_1^2} -1\right)\cr
&& \qquad\,= H^{(2)} + z\,(p_1^2\,q_2^2 - p_2^2\,q_1^2) + o[z^2],
\label{ggg}
\eea  
the nature of the interaction introduced by the non-standard deformation
with respect to (\ref{fg}) can be appreciated. Note that the series
expansion (\ref{ggg}) can be meaningful when the deformation
parameter $z$ is small, and it should be explored in order to analyse the
dynamics.

The $N$-dimensional generalization for this system can be derived 
from the $m$-th order coproducts (\ref{fl}) induced from (\ref{ga}):
\bea  
&&  f_-^{(m)}= \sum_{i=1}^m q_i^2 ,\cr
&&  f_+^{(m)}=\sum_{i=1}^m
\frac {\sinh z q_i^2}{z q_i^2}  p_i^2  e^{z \kk_i^{(m)}(q^2) }, \label{gi}\\
&& f_3^{(m)}=\sum_{i=1}^m
\frac {\sinh z q_i^2}{z q_i}  p_i  e^{z \kk_i^{(m)}(q^2) }  ,
\nonumber
\eea
where the $K$-functions are defined as:
\bea
 \kk_i^{(m)}(x)& =&  - \sum_{k=1}^{i-1}  x_k+ 
\sum_{l=i+1}^m   x_l ,
\label{zzd}\\
 \kk_{ij}^{(m)}(x) & =&\kk_i^{(m)}(x)+ \kk_j^{(m)}(x)\cr
 & =& -(  x_i -   x_j)  - 2\sum_{k=1}^{i-1}   x_k+
2\sum_{l=j+1}^m   x_l,  \qquad i<j .
\label{zze} 
\eea

Therefore, the $N$-dimensional Hamiltonian is just 
$ H^{(N)}_z =  f_+^{(N)}+\para f_-^{(N)}$
and the following constants of motion are deduced:
\be 
  C^{(m)}=- \sum_{i<j}^m 
\frac {\sinh z q_i^2 \sinh z q_j^2}{z^2 q_i^2 q_j^2}
\left({q_i}{p_j} - {q_j}{p_i}\right)^2  e^{ z
\kk_{ij}^{(m)}(q^2)} .
\label{gn}
\ee 
Throughout all the computations leading to (\ref{gn}), the following
property becomes useful:
 \be
\frac {\sinh(z\sum_{i=1}^m   x_i)}z = \sum_{i=1}^m  
\frac {\sinh z   x_i}z  e^{ z \kk_{i}^{(m)}(x)} .
\label{zzf}
\ee
Note that under the limit $z\to 0$ we recover all the expressions
presented in the previous paragraph.


\subsect{Anharmonic chains and their deformation}

We can now consider a more
general dynamical Hamiltonian ${\cal H}$ of the form
\be
{\cal H}=\jp+ {\cal F}(\jm),
\label{anh}
\ee
where ${\cal F}(\jm)$ is an arbitrary smooth function of $\jm$. 
The formalism summarized in Section 2 
ensures that the corresponding system arising from the $N$-th coproduct of
(\ref{anh}) is completely integrable, with ${\cal H}$ being {\it any}
function of the coalgebra generators. Explicitly, this means that any
$N$-particle Hamiltonian of the form
\be 
 H^{(N)}
=    f_+^{(N)} +{\cal F}(f_-^{(N)}) 
=\sum_{i=1}^N p_i^2 + {\cal F}\left(\sum_{i=1}^N q_i^2\right) ,
\label{anha}
\ee
is completely integrable, and (\ref{fk}) are its constants of motion.
Obviously, the integrability (in fact, superintegrability) of (\ref{anha})
is a well-known result, since (\ref{anha}) is just the Hamiltonian
describing the motion of a particle in a $N$-dimensional Euclidean space
under the action of a central potential. In terms of oscillator chains,
the linear case ${\cal F}(\jm)=\para\,\jm$ leads to the previous harmonic
case, and the quadratic one ${\cal F}(\jm)=\jm^2$ would give us an
interacting chain of quartic oscillators. Further definitions of the
function ${\cal F}$ would give us many other anharmonic chains, all of
them sharing the same dynamical symmetry and the same integrals of the
motion.

Let us now stress that, by using the non-standard quantum $sl(2,\Re)$
algebra, it is clear that a realization of (\ref{anh}) in terms of
(\ref{gi}) gives us: \bea 
&& H_z^{(N)} 
=\sum_{i=1}^N
\frac {\sinh z q_i^2}{z q_i^2}  p_i^2  e^{z \kk_i^{(N)}(q^2) }
 + {\cal F}\left(\sum_{i=1}^N q_i^2\right)\cr
&&\qquad\quad =
\sum_{i=1}^N p_i^2 
 + {\cal F}\left(\sum_{i=1}^N q_i^2\right) + 
\sum_{i=1}^N{p_i^2 \,\left(
\frac {\sinh z q_i^2}{z q_i^2}   e^{z \kk_i^{(N)}(q^2) }-1\right) } ,
\label{anhaz}
\eea
which is an integrable deformation of (\ref{anha}) with (\ref{gn}) being
again the associated integrals. This result can be interpreted as a
perturbation of the original anharmonic chain through long-range
interacting terms depending on the momenta.


\sect{The Gaudin Hamiltonian and $sl(2,\Re)$ coalgebras}

If we substitute the canonical realizations used until now in terms of
angular momentum realizations of the same abstract $sl(2,\Re)$ Poisson
coalgebra, the very same construction will lead us to a long-range
interacting ``classical spin chain" of the Gaudin type on which the
quantum deformation can be easily implemented.   

In particular, let us consider the classical angular momentum realization
$S$
\be
g_3^{(1)}=S(\jj)=\s_3^1 ,\qquad
g_+^{(1)}=S(\jp)=\s_+^1 ,\qquad
g_-^{(1)}=S(\jm)=\s_-^1 ,
\label{ha}
\ee
where the  variables $\sigma_i^1$ fulfil
\be
\{\sigma_3,\sigma_+\}=2\sigma_+ ,\qquad 
\{\sigma_3,\sigma_-\}=-2\sigma_- ,\qquad \{\sigma_-,\sigma_+\}=4\sigma_3 ,
\label{pb}
\ee
and are
constrained by a given constant value of the Casimir function (\ref{hgb})
in the form $c_1=(\s_3^1)^2 - \s_-^1 \s_+^1$. 

As usual, $m$ different copies of
(\ref{ha}) (that, in principle, could have different values $c_i$ of the
Casimir)  can be distinguished with the aid of a superscript $\s_l^i$.
Then, the  $m$-th coproduct provides the following realization
of the non-deformed $sl(2,\Re)$ Poisson coalgebra :
\be
 g_l^{(m)}= (S\otimes \dots^{m)}\otimes
S)(\Delta^{(m)}(\s_l))=\sum_{i=1}^{m}{\s_l^i} ,
\qquad l=+,-,3.
\label{hb}
\ee

Now, we can apply the usual construction and take ${\cal H}$ from
(\ref{hsl}). As a consequence, the uncoupled oscillator chain (\ref{fj})
is equivalent to
\be 
H^{(N)} 
=  g_+^{(N)}+\para g_-^{(N)}=\sum_{i=1}^{m}{\left( \s_+^i + \para \s_-^i
\right)} ,
\label{hc}
\ee
and the Casimirs $C^{(m)}$ read $(m=2,\dots, N)$: 
\be 
 C^{(m)}= (g_3^{(m)})^2 - g_-^{(m)} g_+^{(m)} 
  = \sum_{i=1}^{m}{c_i} + 
 \sum_{i<j}^{m} (\s_3^i \s_3^j- \s_-^i \s_+^j -  \s_-^j \s_+^i ).
\label{hd}
\ee 
Note that these Casimirs are just Gaudin Hamiltonians of the hyperbolic
type \cite{Gau,EEKT} (in fact, we shall consider $C^{(N)}$ as the one
defining a general Gaudin magnet).  

The implementation of a non-standard deformation in the Gaudin system is
now straightforward. The deformed angular momentum realization
corresponding to the non-standard deformation $U_z(sl(2,\Re))$ is:
\bea
&& g_-^{(1)}=S_z(\jm)=\s_-^1 ,\cr
&& g_+^{(1)}=S_z(\jp)=\frac {\sinh z \s_-^1}{z
\s_-^1} \s_+^1 ,\label{ia} \\
&& g_3^{(1)}=S_z(\jj)=\frac {\sinh z
\s_-^1}{z \s_-^1}\s_3^1 ,
\nonumber
\eea
where the classical coordinates $\s_l^1$ are defined on the cone  
$c_1=(\s_3^1)^2 -\s_-^1 \s_+^1=0$ (therefore, we are considering the zero
realization). 

It is easy to check that the $m$-th order of the coproduct (\ref{ga}) in
the above representation leads to the following functions
\bea  
 && g_-^{(m)}= \sum_{i=1}^m \s_-^i,\cr
 && g_+^{(m)}=\sum_{i=1}^m
\frac {\sinh z \s_-^i}{z \s_-^i} \s_+^i  e^{z \kk_i^{(m)}(\s_-) }, 
\label{ib}\\
&&g_3^{(m)}=\sum_{i=1}^m
\frac {\sinh z \s_-^i}{z \s_-^i} \s_3^i  e^{z \kk_i^{(m)}(\s_-) } ,
\nonumber
\eea
that define the non-standard deformation of (\ref{hc}). Therefore, the
following ``non-standard Gaudin Hamiltonians" are obtained
\bea
&& C^{(m)}_z=(g_3^{(m)})^2 - \frac {\sinh z g_-^{(m)}}{z} g_+^{(m)}\cr
&&\qquad\quad=
\sum_{i=1}^m \left(\frac {\sinh z \s_-^i}{z \s_-^i}\right)^2
 e^{2 z \kk_i^{(m)}(\s_-) }\left\{ (\s_3^i)^2 -\s_-^i \s_+^i\right\}\cr
&&\quad\qquad\quad +\sum_{i<j}^m \frac {\sinh z \s_-^i \sinh z \s_-^j}{z^2
\s_-^i \s_-^j}
 e^{  z \kk_{ij}^{(m)}(\s_-) } (\s_3^i \s_3^j - \s_-^i \s_+^j
- \s_+^i \s_-^j) .
\label{ic}
\eea
Since we are working in the zero representation with
$(\s_3^i)^2 -\s_-^i \s_+^i=0$, the expression (\ref{ic}) can be simplified 
\be
 C^{(m)}_z=\sum_{i<j}^m \frac {\sinh z \s_-^i \sinh z \s_-^j}{z^2
\s_-^i
\s_-^j}
 e^{  z \kk_{ij}^{(m)}(\s_-) } (\s_3^i \s_3^j - \s_-^i \s_+^j 
- \s_+^i \s_-^j) .
\label{id}
\ee
These integrals are the angular momentum counterparts to
(\ref{gn}). Note that, as it was found for the standard case in \cite{BR},
the deformation can be interpreted as the introduction of a variable range
exchange \cite{Inoz} in the model (compare (\ref{id}) with (\ref{hd})).



{\section*{Acknowledgments}} 

\noindent
 The authors have been partially supported by DGICYT (Project 
PB94--1115) from the Ministerio de Educaci\'on y Ciencia de Espa\~na and
by Junta de Castilla y Le\'on (Projects CO1/396 and CO2/297) and
gratefully acknowledge the invitation to participate in the ``Third
International Workshop on Classical and Quantum Integrable Systems", and
the warm hospitality during their stay in Yerevan. 

\bigskip


\end{document}